\documentclass[12pt]{amsart}
\usepackage{amsmath,amssymb,graphicx}

\newtheorem{propo}{Proposition}[section]
\newtheorem{corol}[propo]{Corollary}
\newtheorem{theor}[propo]{Theorem}
\newtheorem{lemma}[propo]{Lemma}

\theoremstyle{definition}
\newtheorem{defin}[propo]{Definition}
\newtheorem{examp}[propo]{Example}

\theoremstyle{remark}
\newtheorem{remar}[propo]{Remark}

\numberwithin{equation}{section}

\newcommand{\NN }{\mathbb{N}}

\newcommand{\RR }{\mathbb{R}}

\newcommand{\ZZ }{\mathbb{Z}}

\newcommand{\id }{\mathrm{id}}

\newcommand{\Ac }{\mathcal{A}}
\newcommand{\Kc }{\mathcal{K}}

\newcommand{\Cc }{\mathcal{C}}

\newcommand{\rsC }{\mathcal{R}}
\newcommand{\rsCre }{\mathcal{R}^{re}}

\newcommand{\Cm }{C}
\newcommand{\cm }{c}
\newcommand{\rfl }{\rho }
\newcommand{\s }{\sigma }
\DeclareMathOperator{\Aut}{Aut}
\DeclareMathOperator{\End}{End}
\DeclareMathOperator{\Hom}{Hom}
\newcommand{\Wg }{\mathcal{W}}
\newcommand{\Ob }{\mathrm{Ob}}

\newcommand{\re }{^\mathrm{re}}
\newcommand{\rer }[1]{(R\re)^{#1}}

\newcommand{\Addint}{Crystallographic }
\newcommand{\addint}{crystallographic }
\newcommand{\addintn}{crystallographic}

\title[Crystallographic arrangements]
{Crystallographic arrangements: Weyl groupoids and simplicial arrangements}

\author{M.~Cuntz}
\address{Michael Cuntz,
Fachbereich Mathematik,
Universit\"at Kaiserslau\-tern,
Postfach 3049,
D-67653 Kaiserslautern, Germany}
\email{cuntz@mathematik.uni-kl.de}

\begin{document}

\begin{abstract}
We introduce the simple notion of a ``\addint arrangement'' and
prove a one-to-one correspondence between these arrangements and
the connected simply connected Cartan schemes for which the real roots are a finite
root system (up to equivalence on both sides).
We thus obtain a more accessible definition for this very large
subclass of the class of simplicial arrangements for which
a complete classification is known.
\end{abstract}
\maketitle

\section{Introduction}

Originally, Weyl groupoids were introduced as the underlying symmetry structure
of certain Nichols algebras with the purpose to understand pointed Hopf algebras via
the Lifting method of Andruskiewitsch and Schneider \cite{p-AS-02}.
For example the upper triangular part of a small quantum group,
also called Frobenius-Lusztig kernel, is such a Nichols algebra.
Heckenberger classified these structures in the case of Nichols algebras
of diagonal type \cite{p-H-06}.
Later, this invariant was defined in a very general setting
\cite{p-HS-08}.
It plays a similar role as the Weyl group plays for semisimple Lie algebras
and algebraic groups.

After an axiomatic approach to Weyl groupoids had been initiated
in \cite{p-HY-08} and formulated in terms of categories in \cite{p-CH09a},
Heckenberger and the author started the classification of the finite case
\cite{p-CH09b}, \cite{p-CH09c}.
The case of rank three \cite{p-CH09c} revealed a connection to
geometric combinatorics: Heckenberger and Welker \cite{p-HW-10} proved that
the ``root systems'' obtained from Weyl groupoids define
a simplicial arrangement in the same way as reflection groups do.
For example in rank three, 53 of the 67 sporadic simplicial arrangements in
the large component of the Hasse diagram in \cite{p-G-09} may be obtained
from Weyl groupoids (notice that the Weyl {\it groups} only yield $2$ arrangements,
type $A_3$ and $B_3$).
The classification of simplicial arrangements in the
real projective plane is still an open problem and yet there is no other large subclass
of arrangements for which a complete classification is known. This considerably
increases the importance of the theory of Weyl groupoids which henceforth
may be considered as a subject in its own right.

However, the axioms in \cite{p-HY-08} or \cite{p-CH09a} are rather complicated
and seem to be somewhat artificial.
The aim of this article is to give an equivalent, more intuitive definition of ``connected
simply connected Cartan schemes for which the real roots form a finite root
system'' in terms of simplicial arrangements. It turns out that almost two pages
of definitions and axioms (see Section~\ref{gen_res}) reduce to one single axiom (see (I) below).
This shows that Weyl groupoids are a very natural structure after all.
Further, it is conceivable that this new definition could help to find a good
generalization and finally to obtain a classification of the complete class of
simplicial arrangements or at least of the class of free arrangements.

We briefly sketch our result.

Let $\Ac$ be a simplicial arrangement in $\RR^r$, i.e.~the complement of
the union of all hyperplanes in $\Ac$ decomposes into open simplicial cones
called the chambers.
For each hyperplane of $\Ac$, choose a normal vector. Let $R$
be the set of all these normal vectors including their negatives.
For each chamber\footnote{We use the letter $K$ because we will need the
$C$ for ``Cartan matrices'' later.} $K$ of $\Ac$ we have a basis $B^K$ of $\RR^r$
consisting of the normal vectors (in $R$) of the walls of $K$ pointing to the inside
(see Section~\ref{addint_arrs} for more details).
Assume that
\begin{itemize}
\item[(I)] \quad $R \subseteq \sum_{\alpha \in B^K} \ZZ \alpha$\quad for all chambers $K$.
\end{itemize}
We call a pair $(\Ac,R)$ with property (I)
a {\it \addint arrangement}\footnote{Since the word ``integral'' is already in use
for arrangements, ``\addintn'' seems to be the most appropriate.}.
The main result of this article is (see Theorem~\ref{main_thm}):
\begin{theor}
There is a one-to-one correspondence between \addint arrangements 
and connected simply connected Cartan schemes for which the
real roots are a finite root system (up to equivalence on both sides):
\begin{enumerate}
\item Let $(\Ac,R)$ be a \addint arrangement.
Then $R$ is the set of real roots at some object of a connected Cartan scheme $\Cc$.
\item Every connected simply connected Cartan scheme for which the real roots
are a finite root system comes from a unique \addint arrangement up to equivalence.
\end{enumerate}
\end{theor}

In particular, the classification of ``finite Weyl groupoids'' \cite{p-CH10} applies to
our new structure and we obtain a complete classification of \addint arrangements
in arbitrary dimension.
For example, it includes the crystallographic root systems.

This paper is organized as follows. We start with a section in which we
define the \addint arrangements and prove the key lemma (Lemma~\ref{sigma1}).
In the following section we recall the definition of Cartan schemes, Weyl groupoids
and root systems following \cite{p-CH09a}.
In Section~\ref{const_grpd} we construct a Cartan scheme and a root system
for a given \addint arrangement.
In the last section we summarize the results.

\textbf{Acknowledgement.} I am very thankful to I.~Heckenberger and to G.~Malle
for many valuable comments on a previous version of the manuscript.

\section{\Addint and additive arrangements}\label{addint_arrs}

Let $r\in\NN$, $V:=\RR^r$. For $\alpha\in V^*$ we write
$\alpha^\perp = \ker(\alpha)$.
We first recall the definition of a simplicial arrangement (compare \cite[1.2, 5.1]{OT}).
\begin{defin}\label{A_R}
Let $\Ac$ be a simplicial arrangement in $V$, i.e.~
$\Ac=\{H_1,\ldots,H_n\}$ where $H_1,\ldots,H_n$ are distinct linear hyperplanes in $V$
and every component of $V\backslash \bigcup_{H\in\Ac} H$ is an open simplicial cone
(a cone as in equation (\ref{simp_cone}) for a basis
$\alpha^\vee_1,\ldots,\alpha^\vee_r$ of $V$).
Let $\Kc(\Ac)$ be the set of connected components of $V\backslash \bigcup_{H\in\Ac} H$;
they are called the {\it chambers} of $\Ac$.

For each $H_i$, $i=1,\ldots,n$ we choose an element $x_i\in V^*$ such that $H_i=x_i^\perp$.
Let
\[ R = \{\pm x_1,\ldots,\pm x_n\}\subseteq V^*. \]
Further, for $\alpha\in R$ we write
\[ N(\alpha) = \{ v\in V \mid \alpha(v)\ge 0 \}. \]
For each chamber $K\in\Kc(\Ac)$ set
\begin{eqnarray*}
W^K &=& \{ H\in\Ac \mid \dim(H\cap \overline K) = r-1 \}, \\
B^K &=& \{ \alpha\in R \mid \alpha^\perp\in W^K,\quad
N(\alpha) \cap K = K \} \subseteq R.
\end{eqnarray*}
Here, $\overline{K}$ denotes the closure of $K$.
The elements of $W^K$ are the {\it walls} of $K$ and
$B^K$ ``is'' the set of normal vectors of the walls of $K$ pointing to the inside.
Note that
\[ \overline K = \bigcap_{\alpha\in B^K} N(\alpha). \]
Moreover, if $\alpha^\vee_1,\ldots,\alpha^\vee_r$ is the dual basis to
$B^K=\{\alpha_1,\ldots,\alpha_r\}$, then
\begin{equation}\label{simp_cone}
K = \Big\{ \sum_{i=1}^r a_i\alpha_i^\vee \mid a_i> 0 \quad\mbox{for all}\quad
i=1,\ldots,r \Big\}.
\end{equation}
\end{defin}

The following simple lemma is fundamental for the theory.

\begin{lemma}\label{pos_cone}
Let $K\in \Kc(\Ac)$.
If $B^K=\{\alpha_1,\ldots,\alpha_r\}$, then for all $\alpha$ in $R$ we have
$\alpha \in \pm \sum_i \RR_{\ge 0} \alpha_i$.
\end{lemma}

\begin{proof}
Now let $\alpha = \sum_i a_i\alpha_i \in R$ for some $a_1,\ldots,a_r\in\RR$.
Write $\{1,\ldots,r\}=I_1\cup I_2\cup I_3$ such that
$a_i>0$ for $i\in I_1$, $a_i<0$ for $i\in I_2$ and $a_i=0$ for $i\in I_3$.
Assume that $I_1\ne\emptyset\ne I_2$.
Let $\alpha^\vee_1,\ldots,\alpha^\vee_r\in V$ be the dual basis to
$\alpha_1,\ldots,\alpha_r$ and define
\[ v := \sum_{i\in I_1} \frac{1}{|I_1|a_i}\alpha^\vee_i
-\sum_{i\in I_2} \frac{1}{|I_2|a_i}\alpha^\vee_i +
\sum_{i\in I_3} \alpha^\vee_i. \]
Then $\alpha(v)=0$, thus $v\in \alpha^\perp$.
But for all $\beta\in B^K$, $\beta(v)>0$, hence $v\in K\cap \alpha^\perp$
which is a contradiction to the fact that $K$ is a chamber and $\alpha^\perp$
a hyperplane of $\Ac$.
\end{proof}

\begin{defin}
Let $\Ac$ be a simplicial arrangement and $R\subseteq V^*$ a finite set
such that $\Ac = \{ \alpha^\perp \mid \alpha \in R\}$ and $\RR\alpha\cap R=\{\pm \alpha\}$
for all $\alpha \in R$.
For $K\in\Kc(\Ac)$ set
\[ R^K_+ = R \cap \sum_{\alpha \in B^K} \RR_{\ge 0} \alpha. \]
We call $(\Ac,R)$ a {\it \addint arrangement} if
for all $K\in\Kc(\Ac)$:
\begin{itemize}
\item[(I)] \quad $R \subseteq \sum_{\alpha \in B^K} \ZZ \alpha$.
\end{itemize}
We call $(\Ac,R)$ an {\it additive arrangement} if
for all $K\in\Kc(\Ac)$:
\begin{itemize}
\item[(A)] \quad For all $\alpha\in R^K_+\backslash B^K$ there exist
$\beta,\gamma\in R^K_+$ with $\alpha=\beta+\gamma$.
\end{itemize}
Two crystallographic arrangements $(\Ac,R)$, $(\Ac',R')$ in $V$ are called {\it equivalent}
if there exists $\psi\in\Aut(V^*)$ with $\psi(R)=R'$.
We write $(\Ac,R)\cong(\Ac',R')$.
\end{defin}

\begin{remar}\label{I_pos}
By Lemma~\ref{pos_cone}, axiom (I) implies
$R \subseteq \pm \sum_{\alpha \in B^K} \NN_0 \alpha$ for all chambers $K$.
\end{remar}

\begin{examp}\label{coscorf}
Let $\Cc$ be a connected Cartan scheme and assume that $\rsCre(\Cc)$
is a finite root system (see Section~\ref{gen_res}). Let $a$ be an object of $\Cc$.
Then $(\Ac=\{\alpha^{\perp} \mid \alpha\in R_+^a \},R^a)$ is an additive,
\addint arrangement, where we embed $R_+^a$ into the dual space $V^*\cong \RR^r$.
We will restate this in Section~\ref{corresp}, it is a result
of \cite{p-HW-10} and \cite{p-CH09c}.
\end{examp}

\begin{examp}
Crystallographic reflection arrangements (arrangements from Weyl groups) are a special
case of example \ref{coscorf}.
Note that as crystallographic arrangements, arrangements of type $B_r$ and $C_r$ are not
equivalent in rank $r \ge 3$. The reader may want to compare the subsequent development
with the theory of cristallographic reflection arrangements.
\end{examp}

\begin{lemma}\label{trivlem}
Let $(\Ac,R)$ be as above.
Then for all $K\in\Kc(\Ac)$:
\begin{enumerate}
\item $R=R^K_+\dot\cup-R^K_+$,
\item $R^K_+\cap\RR\alpha=\{\alpha\}$ for all $\alpha\in R^K_+$.
\end{enumerate}
Further, if $(\Ac,R)$ is additive then it is \addintn.
\end{lemma}

\begin{proof}
This follows immediately from the definition.
\end{proof}

Our goal is now to prove that any \addint arrangement comes from a ``finite Weyl grou\-poid'',
and that any ``finite Weyl grou\-poid'' is given by some \addint arrangement.

In the sequel, let $(\Ac,R)$ be a \addint arrangement and $K_0,K$ be adjacent chambers,
i.e.~$W^K\cap W^{K_0}=\{\alpha_1^\perp\}$ for some $\alpha_1\in R^{K_0}_+$.
Let $B^{K_0}=\{\alpha_1,\ldots,\alpha_r\}$.

\begin{lemma}\label{pos_sum}
Let $\beta\in B^K$. Then
\[ \beta=-\alpha_1 \quad\mbox{or}\quad \beta \in \sum_{i=1}^r\NN_0 \alpha_i.\]
\end{lemma}

\begin{proof}
Let $\beta\in B^K$ and assume $\beta\in-\sum_{i=1}^r\NN_0 \alpha_i$ (using Rem.~\ref{I_pos}).
Since $N(\gamma)\cap N(\delta)\subseteq N(\gamma+\delta)$ for all
$\gamma,\delta\in V^*$, we have
$\overline K_0=\bigcap_{i=1}^r N(\alpha_i)\subseteq N(-\beta)$.
But $\overline K=\bigcap_{\gamma\in B^K} N(\gamma)$ and thus we obtain
\[ \overline K_0\cap \overline K\subseteq N(\beta)\cap N(-\beta)=\beta^\perp. \]
Since $\alpha_1^\perp$ is a wall of $K_0$ and $K$, we have
$\dim(\overline K_0\cap \overline K)=r-1$ and $\overline K_0\cap \overline K\subseteq\alpha_1^\perp$.
We conclude $\alpha_1^\perp=\beta^\perp$.
\end{proof}

The following lemma will be the key for all arguments concerning
\addint arrangements in this paper.

\begin{lemma}\label{sigma1}
Write $B^K=\{-\alpha_1,\beta_2,\ldots,\beta_r\}$ for some $\beta_2,\ldots,\beta_r\in R$.
Then there exists a permutation $\tau\in S_r$ with $\tau(1)=1$ and such that
\[ \beta_i = c_{\tau(i)}\alpha_1+\alpha_{\tau(i)}, \quad \quad i=2,\ldots,r \]
for certain $c_2,\ldots,c_r\in \NN_0$.
\end{lemma}
\begin{proof}
Let $\sigma$ be the linear map
\[ \sigma : V\rightarrow V,\quad \alpha_1 \mapsto -\alpha_1,\quad
\alpha_i\mapsto\beta_i \quad \text{for} \quad i=2,\ldots,r. \]
With respect to the basis $B^{K_0}$, $\sigma$ is a matrix of the form
\[ \begin{pmatrix}
-1 & c_2 & \cdots & c_r \\
0  &     &        &     \\
\vdots & & A & \\
0 & & &
\end{pmatrix} \]
for some $c_2,\ldots,c_r\in\NN_0$ and $A\in \NN_0^{(r-1)\times (r-1)}$
by Lemma~\ref{pos_sum}.
But by the same argument with the role of $K_0$ and $K$ interchanged, with respect to
$\{-\alpha_1,\beta_2,\ldots,\beta_r \}$,
the $\alpha_2,\ldots,\alpha_r$ have coefficients in $\NN_0$, thus
$A^{-1}\in\NN_0^{(r-1)\times(r-1)}$.
This is clearly only possible if $A$ is a permutation matrix.
\end{proof}

\begin{remar}
Lemma~\ref{sigma1} may also be proved for an arrangement $(\Ac,R)$ for which
a slightly stronger version of axiom (A) holds only at the chamber $K_0$.
But it is not clear whether this stronger axiom transports from $K_0$ to
adjacent chambers, so we do not elaborate on this.
\end{remar}

\begin{defin}
Let $K,K'$ be adjacent chambers and $\{\alpha^\perp\}=W^K\cap W^{K'}$ for
$\alpha\in B^K$. Then by Lemma~\ref{sigma1} there exist
unique $c_\beta\in\NN_0$, $\beta\in B^K\backslash\{\alpha\}$ such that
\[ \varphi_{K,K'} : B^K\rightarrow B^{K'}, \quad
\beta \mapsto c_\beta \alpha + \beta, \quad
\alpha \mapsto -\alpha \]
is a bijection.
\end{defin}

\begin{defin}\label{seq_chain}
Now fix a chamber $K_0$ and an ordering $B^{K_0}=\{\alpha_1,\ldots,\alpha_r\}$.
Then for any sequence $\mu_1,\ldots,\mu_m\in\{1,\ldots,r\}$ we get a unique chain
of chambers $K_0,\ldots,K_m$ such that $K_{i-1}$ and $K_i$ are adjacent and
\[ W^{K_{i-1}}\cap W^{K_i} =
\{\varphi_{K_{i-1},K_i}\ldots\varphi_{K_0,K_1}(\alpha_{\mu_i})^\perp \} \]
for $i=1,\ldots,m$. We will write 
\[ \alpha^{K_i}_{j} := \varphi_{K_{i-1},K_i}\ldots\varphi_{K_0,K_1}(\alpha_{j}) \]
for $i=1,\ldots,m$ and $j=1,\ldots,r$. For $i=0$ we set $\alpha^{K_0}_{j}:=\alpha_{j}$.
But notice that this depends on the chosen chain from $K_0$ to $K_i$, that
is, the $\mu_1,\ldots,\mu_i$.

For $i=0,\ldots,m-1$, denote by $\sigma^{K_i}_{\mu_{i+1}}$ the linear map given by
\[ \sigma^{K_i}_{\mu_{i+1}}(\alpha) = \varphi_{K_i,K_{i+1}}(\alpha) \]
for all $\alpha\in B^{K_i}$. By Lemma~\ref{sigma1} (we have $A=\id$ now),
\[ \sigma^{K_i}_{\mu_{i+1}}\in\Aut(V^*), \quad {(\sigma^{K_i}_{\mu_{i+1}}})^2=\id. \]
Note that $\sigma^{K_i}_{\mu_{i+1}}$ is a reflection.
For a sequence $\mu_1,\ldots,\mu_m$ we may abbreviate
\[ \sigma_{\mu_m}\ldots\sigma_{\mu_2} \sigma^{K_0}_{\mu_1}
:= \sigma^{K_{m-1}}_{\mu_m}\ldots\sigma^{K_1}_{\mu_2} \sigma^{K_0}_{\mu_1} \]
since $K_1,\ldots,K_m$ are uniquely determined by the sequence of $\mu_i$.
\end{defin}

\section{Cartan schemes and Weyl groupoids}\label{gen_res}

We briefly recall the notions of Cartan schemes, Weyl groupoids
and root systems, following \cite{p-CH09a,p-CH09b}.
The foundations of the general theory have been developed in
\cite{p-HY-08} using a somewhat different terminology.

\begin{defin}
Let $I$ be a non-empty finite set and
$\{\alpha_i\,|\,i\in I\}$ the standard basis of $\ZZ ^I$.
By \cite[\S 1.1]{b-Kac90} a {\it generalized Cartan matrix}
$\Cm =(\cm _{ij})_{i,j\in I}$
is a matrix in $\ZZ ^{I\times I}$ such that
\begin{enumerate}
  \item[(M1)] $\cm _{ii}=2$ and $\cm _{jk}\le 0$ for all $i,j,k\in I$ with
    $j\not=k$,
  \item[(M2)] if $i,j\in I$ and $\cm _{ij}=0$, then $\cm _{ji}=0$.
\end{enumerate}
\end{defin}

\begin{defin}
Let $A$ be a non-empty set, $\rfl _i : A \to A$ a map for all $i\in I$,
and $\Cm ^a=(\cm ^a_{jk})_{j,k \in I}$ a generalized Cartan matrix
in $\ZZ ^{I \times I}$ for all $a\in A$. The quadruple
\[ \Cc = \Cc (I,A,(\rfl _i)_{i \in I}, (\Cm ^a)_{a \in A})\]
is called a \textit{Cartan scheme} if
\begin{enumerate}
\item[(C1)] $\rfl _i^2 = \id$ for all $i \in I$,
\item[(C2)] $\cm ^a_{ij} = \cm ^{\rfl _i(a)}_{ij}$ for all $a\in A$ and
  $i,j\in I$.
\end{enumerate}
\end{defin}

\begin{defin}
  Let $\Cc = \Cc (I,A,(\rfl _i)_{i \in I}, (\Cm ^a)_{a \in A})$ be a
  Cartan scheme. For all $i \in I$ and $a \in A$ define $\s _i^a \in
  \Aut(\ZZ ^I)$ by
  \begin{align}
    \s _i^a (\alpha_j) = \alpha_j - \cm _{ij}^a \alpha_i \qquad
    \text{for all $j \in I$.}
    \label{eq:sia}
  \end{align}
  The \textit{Weyl groupoid of} $\Cc $
  is the category $\Wg (\Cc )$ such that $\Ob (\Wg (\Cc ))=A$ and
  the morphisms are compositions of maps
  $\s _i^a$ with $i\in I$ and $a\in A$,
  where $\s _i^a$ is considered as an element in $\Hom (a,\rfl _i(a))$.
  The category $\Wg (\Cc )$ is a groupoid in the sense that all morphisms are
  isomorphisms.
\end{defin}
  As above, for notational convenience we will often neglect upper indices referring to
  elements of $A$ if they are uniquely determined by the context. For example,
  the morphism $\s _{i_1}^{\rfl _{i_2}\cdots \rfl _{i_k}(a)}
  \cdots \s_{i_{k-1}}^{\rfl _{i_k(a)}}\s _{i_k}^a\in \Hom (a,b)$,
  where $k\in \NN $, $i_1,\dots,i_k\in I$, and
  $b=\rfl _{i_1}\cdots \rfl _{i_k}(a)$,
  will be denoted by $\s _{i_1}\cdots \s _{i_k}^a$ or by
  $\id _b\s _{i_1}\cdots \s_{i_k}$.
  The cardinality of $I$ is termed the \textit{rank of} $\Wg (\Cc )$.
\begin{defin}
  A Cartan scheme is called \textit{connected} if its Weyl grou\-poid
  is connected, that is, if for all $a,b\in A$ there exists $w\in \Hom (a,b)$.
  The Cartan scheme is called \textit{simply connected},
  if $\Hom (a,a)=\{\id _a\}$ for all $a\in A$.

  Two Cartan schemes $\Cc =\Cc (I,A,(\rfl _i)_{i\in I},(\Cm ^a)_{a\in A})$
  and $\Cc '=\Cc '(I',A',$
  $(\rfl '_i)_{i\in I'},({\Cm '}^a)_{a\in A'})$
  are termed \textit{equivalent}, if there are bijections $\varphi _0:I\to I'$
  and $\varphi _1:A\to A'$ such that
  \begin{align}\label{eq:equivCS}
    \varphi _1(\rfl _i(a))=\rfl '_{\varphi _0(i)}(\varphi _1(a)),
    \qquad
    \cm ^{\varphi _1(a)}_{\varphi _0(i) \varphi _0(j)}=\cm ^a_{i j}
  \end{align}
  for all $i,j\in I$ and $a\in A$. We write then $\Cc\cong \Cc'$.
\end{defin}

  Let $\Cc $ be a Cartan scheme. For all $a\in A$ let
  \[ \rer a=\{ \id _a \s _{i_1}\cdots \s_{i_k}(\alpha_j)\,|\,
  k\in \NN _0,\,i_1,\dots,i_k,j\in I\}\subseteq \ZZ ^I.\]
  The elements of the set $\rer a$ are called \textit{real roots} (at $a$).
  The pair $(\Cc ,(\rer a)_{a\in A})$ is denoted by $\rsC \re (\Cc )$.
  A real root $\alpha\in \rer a$, where $a\in A$, is called positive
  (resp.\ negative) if $\alpha\in \NN _0^I$ (resp.\ $\alpha\in -\NN _0^I$).
  In contrast to real roots associated to a single generalized Cartan matrix,
  $\rer a$ may contain elements which are neither positive nor negative. A good
  general theory, which is relevant for example for the study of Nichols
  algebras, can be obtained if $\rer a$ satisfies additional properties.

\begin{defin}
  Let $\Cc =\Cc (I,A,(\rfl _i)_{i\in I},(\Cm ^a)_{a\in A})$ be a Cartan
  scheme. For all $a\in A$ let $R^a\subseteq \ZZ ^I$, and define
  $m_{i,j}^a= |R^a \cap (\NN_0 \alpha_i + \NN_0 \alpha_j)|$ for all $i,j\in
  I$ and $a\in A$. We say that
  \[ \rsC = \rsC (\Cc , (R^a)_{a\in A}) \]
  is a \textit{root system of type} $\Cc $, if it satisfies the following
  axioms.
  \begin{enumerate}
    \item[(R1)]
      $R^a=R^a_+\cup - R^a_+$, where $R^a_+=R^a\cap \NN_0^I$, for all
      $a\in A$.
    \item[(R2)]
      $R^a\cap \ZZ\alpha_i=\{\alpha_i,-\alpha_i\}$ for all $i\in I$, $a\in A$.
    \item[(R3)]
      $\s _i^a(R^a) = R^{\rfl _i(a)}$ for all $i\in I$, $a\in A$.
    \item[(R4)]
      If $i,j\in I$ and $a\in A$ such that $i\not=j$ and $m_{i,j}^a$ is
      finite, then
      $(\rfl _i\rfl _j)^{m_{i,j}^a}(a)=a$.
  \end{enumerate}
\end{defin}

  The axioms (R2) and (R3) are always fulfilled for $\rsC \re $.
  The root system $\rsC $ is called \textit{finite} if for all $a\in A$ the
  set $R^a$ is finite. By \cite[Prop.\,2.12]{p-CH09a},
  if $\rsC $ is a finite root system
  of type $\Cc $, then $\rsC =\rsC \re $, and hence $\rsC \re $ is a root
  system of type $\Cc $ in that case.

\section{Constructing a Weyl groupoid and a root system}\label{const_grpd}

We now construct a Weyl groupoid for a given \addint arrangement
$(\Ac,R)$ of rank $r$. Let $I:=\{1,\ldots,r\}$.
Fix a chamber $K_0\in\Kc(\Ac)$ and an ordering $B^{K_0}=\{\alpha_1,\ldots,\alpha_r\}$.
Consider the set
\[ \hat A := \{ (\mu_1,\ldots,\mu_m)
\mid m\in\NN,\:\: \mu_1,\ldots,\mu_m\in I \} \]
and write $a.\nu$ for the sequence $(\mu_1,\ldots,\mu_m,\nu)$ if
$a=(\mu_1,\ldots,\mu_m)$, $\nu\in I$.
We have a map
\[ \pi : \hat A \rightarrow \End(V^*),\quad
(\mu_1,\ldots,\mu_m)\mapsto \sigma_{\mu_m}\ldots\sigma_{\mu_2} \sigma^{K_0}_{\mu_1} \]
which yields an equivalence relation $\sim$ on $\hat A$ via
\[ v\sim w \quad :\Longleftrightarrow \quad \pi(v)=\pi(w) \]
for $v,w\in \hat A$. Let
\[ A := \hat A_{/\sim}. \]
Each sequence $a=(\mu_1,\ldots,\mu_m)\in\hat A$ determines a unique map $\varphi_a$ by
\[ \varphi_a := \varphi_{K_{m-1},K_m}\ldots\varphi_{K_0,K_1} \]
where $K_0,\ldots,K_m$ is the sequence of chambers corresponding to $a$.
\begin{lemma}
The set $A$ is finite.
\end{lemma}
\begin{proof}
Since each sequence $a\in \hat A$ leads to a unique chamber $K$ and
$\varphi_a$ is a bijection $B^{K_0}\rightarrow B^K$, the set $A$ has at
most $|\Kc(\Ac)| r!$ elements and thus is finite.
\end{proof}

\begin{remar}\label{def_Ka}
Two equivalent sequences certainly lead to the same chamber.
For $a\in A$, we will write $K_a$ for this unique chamber.
\end{remar}

\begin{defin}
For each $a\in \hat A$ we construct a Matrix
$C^a\in\ZZ^{r\times r}$ in the following way:

Let $i,j\in I$ and let $K'$ be the chamber adjacent
to $K_a$ with $W^{K_a}\cap W^{K'}=\{\varphi_a(\alpha_i)^\perp\}$.
By Lemma~\ref{sigma1} there exist integers $c_{i,1},\ldots,c_{i,r}$ such that
\[ \varphi_{a.i}(\alpha_j) = -c_{i,j} \varphi_a(\alpha_i)+\varphi_a(\alpha_j). \]
We set
\[ C^a := (c_{i,j})_{1\le i,j \le r}. \]
\end{defin}

\begin{lemma}
Let $a\in \hat A$. Then $C^a$ satisfies {\rm (M1)} and {\rm (M2)}.
\end{lemma}
\begin{proof}
By Lemma~\ref{sigma1}, $c_{i,i}=2$ and $c_{i,j}\le 0$ for $i\ne j$. Thus
the matrix $C^a$ satisfies (M1).

Without loss of generality take $a=()$, $m=0$ and $K_0=K$.
Assume that $c_{i,j}=0$ and let $K',K''$ be the chambers with
$W^K\cap W^{K'}=\{\alpha_i^\perp\}$, $W^K\cap W^{K''}=\{\alpha_j^\perp\}$.
Thus $\sigma^K_i(\alpha_j)=\alpha_j$. But then
\[ R\cap\langle\alpha_i,\alpha_j\rangle = \pm\{\alpha_i,\alpha_j\} \]
(otherwise we would have a hyperplane intersecting a chamber)
and therefore $\sigma^K_j(\alpha_i)\in \pm\{\alpha_i,\alpha_j\}$ which
is only possible if $\sigma^K_j(\alpha_i)=\alpha_i$. Hence $c_{j,i}=0$
and axiom (M2) holds for $C^a$.
\end{proof}
Notice that $C^a=C^b$ for $a,b\in\hat A$ with $\pi(a)=\pi(b)$, so we
obtain a unique Cartan matrix for each element of $A$, which we will
also denote $C^a$, $a\in A$.

We now define the maps $\rho_i : A\rightarrow A$, $i\in I$.
First, set
\[ \hat\rho_i : \hat A \rightarrow \hat A,\quad
a \mapsto a.i. \]
Since $\pi(a)=\pi(b)$ implies $\varphi_a=\varphi_b$ for $a,b\in A$,
this induces well-defined maps
\[ \rho_i : A \rightarrow A, \quad \overline{a} \mapsto \overline{a.i}, \]
where $\overline{a}$ is the equivalence class of $a$ in $A$.
By Lemma~\ref{sigma1},
\[ \rho_i^2=\id, \quad c_{i,j}^a=c_{i,j}^{\rho_i(a)} \]
for all $i,j\in I$, $a\in A$, hence (C1) and (C2) are satisfied.
We obtain:

\begin{propo}\label{CW_A}
Let $(\Ac,R)$ be a \addint arrangement of rank $r$ and $K_0$ a chamber.
Then $I$, $A$, $(\rfl _i)_{i \in I}$, $(\Cm ^a)_{a \in A}$ defined as above
define a Cartan scheme which we denote $\Cc=\Cc(\Ac,R,K_0)$.
Further, we obtain a Weyl groupoid $\Wg(\Cc)=\Wg(\Ac,R,K_0)$.
\end{propo}

\begin{remar}
The Cartan scheme $\Cc=\Cc(\Ac,R,K_0)$ depends on the ordering chosen
on $B^{K_0}$, but a different ordering clearly yields an equivalent Cartan scheme.
\end{remar}

We may now construct a root system of type $\Cc(\Ac,R,K_0)$.
Let $a=(\mu_1,\ldots,\mu_m)\in \hat A$ and let
\[ \phi_a : V^* \rightarrow \RR^r \]
be the coordinate map for elements of $V^*$ with respect to the basis
$\varphi_a(\alpha_1),\ldots,\varphi_a(\alpha_r)$ (in this ordering).
Set
\[ R^a := \phi_a(R). \]
Again, notice that $R^a=R^b$ for $a,b\in \hat A$ with $\pi(a)=\pi(b)$,
so we get a unique set which we also denote $R^a$ for each element of $a\in A$.

\begin{lemma}\label{rank_2_homot}
Let $i,j\in I$, $i\ne j$. Starting with $K_0$,
let $K_0,K_1,\ldots$ be the unique chain given by
$i,j,i,j,\ldots$ Then $K_m=K_0$ for some smallest $m\in\NN$ and
\[ \varphi_{K_{m-1},K_m}\ldots\varphi_{K_1,K_2}\varphi_{K_0,K_1}
= \id. \]
Further,
\[ m = |R\cap \langle\alpha_i,\alpha_j\rangle| \]
where $B^{K_0}=\{\alpha_1,\ldots,\alpha_r\}$.
\end{lemma}
\begin{proof}
By Lemma~\ref{sigma1} and Def.~\ref{seq_chain} the chain $(i,j,i,j,\ldots)$ determines a
sequence of endomorphisms with determinant $-1$ leaving the subspace
$U:=\langle\alpha_i,\alpha_j\rangle$ invariant.

Since we have only finitely many chambers, $K_k=K_m$ for some $k<m\in\NN_0$, $m$ minimal.
Assume that $k>0$. Then since ${\alpha^{K_k}_i}^\perp$ is the wall between $K_{k-1}$ and $K_k$ and
${\alpha^{K_k}_j}^\perp$ is the wall between $K_k$ and $K_{k+1}$ (or possibly with $i,j$ exchanged),
the wall between $K_{m-1}$ and $K_m$ is neither ${\alpha^{K_k}_j}^\perp$ nor ${\alpha^{K_k}_i}^\perp$
unless $K_{m-1}=K_{k-1}$ or $K_{m-1}=K_{k+1}$, contradicting the minimality of $m$.
Thus $K_0=K_m$, $m\in\NN$ minimal.

Denote $V':=V/U^\perp$ and $\kappa : V\rightarrow V'$ the projection.
The set $R\cap U$ yields a simplicial arrangement in $V'$ of rank $2$:
$\{ \kappa(\alpha^\perp) \mid \alpha \in R\cap U \}$.
The chain of chambers $K_0,\ldots,K_m$ gives a chain in $V'$ of consecutively adjacent chambers
$\kappa(K_0),\ldots,\kappa(K_m)$.
But $K_0=K_m$, thus $m$ has to be even and the determinant of the endomorphism is $1$.
We conclude that $\varphi_{K_{m-1},K_m}\ldots\varphi_{K_1,K_2}\varphi_{K_0,K_1}$
is a power of the transposition $(i,j)$ with determinant $1$, hence it is equal to $\id$.

The last assertion also follows by considering the arrangement of rank $2$ in $V'$:
Assume that $\kappa(K_0)=\kappa(K_\nu)$ for some $\nu \in\{1,\ldots,m\}$.
Note that as above $\nu$ has to be even. By Lemma~\ref{sigma1}, with
respect to the basis $B^{K_0}$, the map $\sigma^{K_{\nu-1}}_j\ldots\sigma_i^{K_0}$
is a matrix (without loss of generality $\{i,j\}=\{1,2\}$)
\[ \begin{pmatrix}
     1 & 0 & * & \cdots & * \\
     0 & 1 & * & \cdots & * \\
     0 & 0 & 1 & & 0 \\
\vdots & \vdots &   & \ddots & \\
     0 & 0 & 0 & & 1
\end{pmatrix}. \]
This implies $K_0\subseteq K_\nu$, thus $K_0=K_\nu$. Since $m$ was chosen minimal, $m=\nu$ and
$m$ is the number of chambers of the arrangement in $V'$, exactly one for each
element of $|R\cap \langle\alpha_i,\alpha_j\rangle|$.
\end{proof}

Thus the groupoid $\Wg(\Ac,R,K_0)$ is a Coxeter groupoid in the sense of \cite{p-HY-08}.

\begin{propo}\label{R_A}
The sets $R^a$, $a\in A$ as above are the real roots of $\Cc=\Cc(\Ac,R,K_0)$ and
form a root system of type $\Cc$.
\end{propo}
\begin{proof}
By axiom (I) we have $R^a\subseteq\ZZ^I$. Axioms (R1), (R2) and (R3)
follow immediately from the fact that $\Ac$ is an arrangement and
the definition of $\sigma_i^K$ for a chamber $K$ and $i\in I$.
Remark that these sets are indeed the real roots of $\Cc(\Ac,R,K_0)$
as defined above.
It remains to check axiom (R4). But this is Lemma~\ref{rank_2_homot}.
\end{proof}

\section{The correspondence}\label{corresp}

In this section we summarize the relation between \addint arrangements
and root systems.

\begin{lemma}\label{conn_A}
Let $(\Ac,R)$ be a \addint arrangement.
Then for every chamber $K_0$ the Cartan scheme $\Cc(\Ac,R,K_0)$ is connected.
Moreover, every chamber $K$ is attained by some sequence.
\end{lemma}
\begin{proof}
Let $U$ be the union of all intersections of distinct hyperplanes of $\Ac$.
Then for each pair of chambers $K,K'$ there exists a chain of
consecutively adjacent chambers $K=K_1,\ldots,K_m=K'$,
since $V\backslash U$ is connected.
Thus every chamber $K$ is attained from $K_0$ by some sequence.
Now if $\overline{a}$, $\overline{b}$ are objects in $\Cc(\Ac,R,K_0)$,
i.e.~$a,b\in \hat A$, and $a^{-1}$ denotes the reversed sequence, then
$b a^{-1}$ gives a morphism from $\overline{a}$ to $\overline{b}$,
so $\Cc(\Ac,R,K_0)$ is connected.
\end{proof}

\begin{lemma}\label{sc_A}
Let $(\Ac,R)$ be a \addint arrangement.
Then for every chamber $K_0$ the Cartan scheme $\Cc(\Ac,R,K_0)$ is simply connected.
\end{lemma}
\begin{proof}
This follows from the definition of the map $\pi$.
\end{proof}

\begin{lemma}\label{equiv}
Let $(\Ac,R)$ be a \addint arrangement and $K,K'$ be chambers.
Then $\Cc(\Ac,R,K)$ and $\Cc(\Ac,R,K')$ are equivalent Cartan schemes.
\end{lemma}
\begin{proof}
Choose orderings $B^K=\{\alpha_1,\ldots,\alpha_r\}$ and
$B^{K'}=\{\alpha'_1,\ldots,\alpha'_r\}$ which yield $\Cc(\Ac,R,K)$ and $\Cc(\Ac,R,K')$.
Then to any chain $a=(\mu_1,\ldots,\mu_m)\in \hat A$ of consecutively adjacent chambers
from $K$ to $K'$ belongs a bijection $\varphi_a : B^K \rightarrow B^{K'}$,
thus an element $\varphi_0\in S_r$ since the set of labels is
$\{1,\ldots,r\}$ on both sides. Up to equivalence of Cartan schemes, we may assume
without loss of generality that $\varphi_0=\id$.

If $A$ and $A'$ denote the objects of $\Cc(\Ac,R,K)$ resp.~$\Cc(\Ac,R,K')$, then
\[ \varphi_1 : A' \rightarrow A, \quad b \mapsto \mu_1\ldots\mu_m.b \]
is bijective, since all paths are invertible.
Clearly, $\varphi_1$ is compatible with the Cartan schemes
($\rho$, $\rho'$ and the Cartan matrices) and thus it gives an
equivalence of Cartan schemes $\Cc(\Ac,R,K)\cong\Cc(\Ac,R,K')$.
\end{proof}

We can now state the main theorem.

\begin{theor}\label{main_thm}
Let $\mathfrak A$ be the set of all \addint arrangements and
$\mathfrak C$ be the set of all connected simply connected Cartan schemes for
which the real roots are a finite root system.
Then the map
\[ \Lambda : \mathfrak{A}/_{\cong} \rightarrow \mathfrak{C}/_{\cong}, \quad
\overline{(\Ac,R)} \mapsto \overline{\Cc(\Ac,R,K)}, \]
where $K$ is any chamber of $\Ac$, is a bijection.
\end{theor}
\begin{proof}
The map $\Lambda$ is well-defined by Prop.~\ref{CW_A}, Prop.~\ref{R_A},
Lemma~\ref{conn_A}, Lemma~\ref{sc_A} and Lemma~\ref{equiv}.
The surjectivity of $\Lambda$ is \cite[Corollary 4.7]{p-HW-10}.
If $\Lambda(\overline{(\Ac,R)})=\Lambda(\overline{(\Ac',R')})$
then the ranks are equal, say $r$, and
we have a bijection $\varphi _0\in S_r$ and a bijection $\varphi _1$ between
the objects compatible with the Cartan schemes.
But $\varphi_0$, $\varphi_1$ define an element $\psi\in\Aut({\RR^r}^*)$
which is an equivalence $(\Ac,R)\cong(\Ac',R')$, thus $\Lambda$ is injective.
\end{proof}

Recall the following theorem which is essential for the
classification of finite Weyl groupoids:
\begin{theor}{\cite[Thm.~2.10]{p-CH09c}}\label{root_is_sum}
  Let $\Cc $ be a Cartan scheme. Assume that $\rsC \re (\Cc )$ is a finite
  root system of type $\Cc $.  Let $a\in A$ and $\alpha \in R^a_+$.
  Then either $\alpha$ is simple, or it is the sum of two positive roots.
\end{theor}

\begin{corol}
Let $(\Ac,R)$ be as in Def.~\ref{A_R}. Then $(\Ac,R)$ is \addint if and only if
it is additive.
\end{corol}

\begin{proof}
This is Thm.~\ref{main_thm}, Lemma~\ref{trivlem} and Thm.~\ref{root_is_sum}.
\end{proof}

\bibliographystyle{amsplain}
\bibliography{refs}

\end{document}